\documentclass[12pt,reqno]{amsart}
\usepackage{amsmath,amsthm,amssymb,amsfonts,amscd}
\usepackage{mathrsfs}
\usepackage{bbm}
\usepackage{bbding}
\usepackage{graphicx,latexsym}
\usepackage[backref=page]{hyperref}
\usepackage{hyperref}
\usepackage{geometry}
\usepackage{color}
\usepackage{xcolor}
\usepackage{picture,epic}
\usepackage{tikz}

\usepackage{enumitem}

\numberwithin{equation}{section}

\theoremstyle{plain}
\newtheorem{theorem}{Theorem}[section]
\newtheorem{lemma}[theorem]{Lemma}

\newtheorem{proposition}[theorem]{Proposition}

\theoremstyle{definition}
\newtheorem{definition}[theorem]{Definition}

\theoremstyle{remark}

\renewcommand{\Re}{\operatorname{Re}}

\newcommand{\GL}{\operatorname{GL}}

\newcommand{\dd}{\mathrm{d}}

\makeatletter
\def\@tocline#1#2#3#4#5#6#7{\relax
  \ifnum #1>\c@tocdepth 
  \else
    \par \addpenalty\@secpenalty\addvspace{#2}%
    \begingroup \hyphenpenalty\@M
    \@ifempty{#4}{%
      \@tempdima\csname r@tocindent\number#1\endcsname\relax
    }{%
      \@tempdima#4\relax
    }%
    \parindent\z@ \leftskip#3\relax \advance\leftskip\@tempdima\relax
    \rightskip\@pnumwidth plus4em \parfillskip-\@pnumwidth
    #5\leavevmode\hskip-\@tempdima
      \ifcase #1
       \or\or \hskip 1em \or \hskip 2em \else \hskip 3em \fi%
      #6\nobreak\relax
    \hfill\hbox to\@pnumwidth{\@tocpagenum{#7}}\par
    \nobreak
    \endgroup
  \fi}
\makeatother

\begin{document}

\title{On the second integral moment of $L$-functions}
\author{Liangxun Li}
\address{Data Science Institute and School of Mathematics \\ Shandong University \\ Jinan \\ Shandong 250100 \\China}
\email{lxli@mail.sdu.edu.cn}

\date{\today}

\begin{abstract}
   Assume that the generalized Ramanujan conjecture holds on the automorphic $L$-function $L(s, \pi)$ on $\GL_d$ over $\mathbb{Q}$ with $d\geq 3$, we can obtain a small log-saving non-trivial bound on the second integral moment of $L(1/2+it, \pi)$. Specifically the bound
   \[
        \int_{T}^{2T}\Big|L\big(\frac{1}{2}+it, \pi\big)\Big |^2 \dd t\ll_{\pi} \frac{T^{\frac{d}{2}}}{\log^{\eta_d}T}
   \]
   holds for a small constant $\eta_d>0$. 
\end{abstract}
\keywords{Automorphic form, moment of $L$-function, multiplicative function}
\subjclass[2020]{11F66, 11N75}

\maketitle

\section{Introduction} \label{sec:Intr}
Let $L(s, \pi)$  be the $L$-function associated $\pi$ which is an irreducible cuspidal automorphic representation of $\GL_d(\mathbb{A})$. A fundamental problem is to establish the asymptotic formulas, for the even integral moments of in $t$-aspect:
\begin{equation}
M_{\pi}^{2k}(T):=\int_{0}^{T}\Big|L\Big(\frac{1}{2}+it, \pi\Big)\Big|^{2k}\dd t,
\end{equation}
where $k\in \mathbb{Z}_{\geq 1}$ and $T\geq 2$ large. 
Understanding the growth of such moments offers deep insight into the distribution of values of 
$L$-functions along the critical line. 
However, deriving asymptotics for $M_{\pi}^{2k}(T)$ remains a formidable challenge. 
Currently, such results are known only for a limited number of families of $L$-functions and for small values of $k.$ 
The difficulty increases significantly with the degree $d$ of the $L$-function and the order $2k$ of the moment.

A classical example is the Riemann zeta function $\zeta(s)$, which corresponds to the trivial representation on $\GL_1(\mathbb{A}).$
In this case $(d=1)$, the asymptotic behaviours for the second and fourth moments are known. Specifically, Hardy and Littlewood \cite{Hardy-Littlewood} proved that
\begin{equation}\label{eqn: 2m_zeta}
    M_{\zeta}^2(T):=\int_{0}^{T}\Big|\zeta(\frac{1}{2}+it)\Big|^{2}\dd t\sim T\log T.
\end{equation}
Ingham \cite{Ingham} established that
\begin{equation}
     M_{\zeta}^4(T):=\int_{0}^{T}\Big|\zeta(\frac{1}{2}+it)\Big|^{4}\dd t\sim \frac{1}{2\pi^2}T\log^4 T.
\end{equation}
The asymptotic formulas for higher moments of $\zeta$ remains open for us. 
It  is conjectured that $ M_{\zeta}^{2k}(T)\sim c_k T\log^{k^2} T$ for some constant $c_k$. Such asymptotic formulas have deep connections to random matrix theory,
see e.g.  \cite{KS00}, \cite{CFKRS}.  To support the conjectures, Ramachandra \cite{Ramachandra} established the lower bound $M_{\zeta}^{2k}(T)\gg_k T\log^{k^2} T.$
Conditionally on Riemann Hypothesis (RH), Soundararajan \cite{Sound09moments} proved the upper bound $M_{\zeta}^{2k}(T)\ll_{k, \varepsilon} T\log^{k^2+\varepsilon} T$ for any $\varepsilon>0.$  Applying the refinement of Harper \cite{Harper13}, $\varepsilon$ can be removed. By assuming some additional conjectures, such as additive divisor conjectures and RH, one can get asymptotic formulas for $M_{\zeta}^{6}(T)$ and $M_{\zeta}^{8}(T)$, see \cite{Ng21}, \cite{NSW}. 
For unconditional results, Heath-Brown \cite{HB} gave the strong upper bound estimate $M_{\zeta}^{12}(T)\ll T^{2}(\log T)^{17}.$
By the interpolation on the estimates of $M_{\zeta}^{12}(T)$ and $M_{\zeta}^{4}(T)$, one can get upper bounds for moments of order between  $4$ and $12$:
\[
M_{\zeta}^{2k}(T)\ll_\varepsilon T^{\frac{k}{4}+\frac{1}{2}+\varepsilon}, \quad 2\leq k\leq 6.
\]
It remains open to break the above exponent, especially for the sixth moment.

For the $\rm GL(2)$ automorphic $L$-function $L(s, f)$ (degree $d=2$) where $f$ is a holomorphic or Maass cusp form, an asymptotic formula for the second moment is known:
\begin{equation}
M_{f}^2(T):=\int_{0}^{T}\Big|L\Big(\frac{1}{2}+it, f\Big)\Big|^{2}\dd t\sim c_f T\log T,
\end{equation}
where $c_f$ is a certain constant depending on $f$. For instance, Good \cite{Go82} proved an asymptotic formula with a power-saving error term $O(T^{\frac{2}{3}+\varepsilon})$ in the holomorphic case. 

For the $L$-function $L(s, \pi)$ of degree $d\geq 3$, it is still open to understand the asymptotic behavior for $M_{\pi}^{2k}(T)$.
Assuming the Generalized  Ramanujan conjecture (GRC), Pi \cite{Pi} establish the sharp lower bound $M_\pi^{2k}(T)\gg_{\pi, k} T(\log T)^{k^2}$. 
Additionally on GRH, Tang and Xiao \cite{TX16} showed the sharp upper bound $M_\pi^{2k}(T)\ll_{\pi, k} T(\log T)^{k^2}$.Unconditionally,
even in the simplest case $k=1$, the sharp upper bound is unknown.  
Heuristically, it is expected that $M_{\pi}^{2k}(T)= T^{1+o_{k,\pi}(1)}$ for any fixed $k$ and $\pi$.
By employing the approximate functional equation and applying the integral large sieve inequality, one obtains the trivial bound:
\begin{equation}\label{eqn: M(T)_trivialbound}
    M_{\pi}^{2k}(T)\ll_{\pi, k, \varepsilon} T^{\frac{dk}{2}+\varepsilon}.
\end{equation}
Under the assumption of the Generalized Lindel\"of Hypothesis (GLH), which asserts
\begin{equation}
    L\Big(\frac{1}{2}+it, \pi\Big)\ll_{\pi, \varepsilon} (1+|t|)^{\varepsilon},
\end{equation}
it follows that for any integer $k\geq 1$,
\begin{equation}
     M_{\pi}^{2k}(T)\ll_{\pi, k, \varepsilon} T^{1+\varepsilon}.
\end{equation}
In general, for $\pi$ of the degree $d\geq 3$,  
improving upon the trivial bound \eqref{eqn: M(T)_trivialbound} with a power saving is a profoundly difficult problem. Progress in this direction would also imply a breakthrough in breaking the convexity bound for  $L(s, \pi)$ in $t$-aspect.
For the degree $3$ $L$-function $L(s, \phi)$,
where $\phi$ is a $\rm GL(3)$ Hecke--Maass cusp form,
the first non-trivial bound for $M_{\phi}^2(T)$ is obtained by Pal \cite{Pal22}. He proved that
$
M_{\phi}^2(T)\ll_{\phi, \varepsilon} T^{\frac{3}{2}-\frac{3}{32}+\varepsilon}.
$
This bound has various applications such as the Rankin--Selberg problem and the strong bound for the first moment of $\GL(3)\times \GL(2)$ $L$-function on $\GL(2)$ spectral aspect.
This result was recently strengthened by Dasgupta, Leung and Young \cite{DLY24} who improved the exponent to $\frac{4}{3}$.

In this paper, we consider the second moment $M_{\pi}^2(T)$ for general $\pi$ (of degree $d\geq 3$). Under the assumption of the Generalized Ramanujan Conjecture (GRC), we establish the following new upper bound.
\begin{theorem}\label{thm: 2M_L}
Let $\mathbb{A}$ be the ring of adel\'es of $\mathbb{Q}$ and $\pi$ be an irreducible cuspidal automorphic representation of $\GL_d(\mathbb{A})$ with its degree $d\geq 3$.
Assume that the generalized Ramanujan conjecture (GRC) holds on $L(s, \pi)$,  then for $T\geq 2$ large, we have
\begin{equation}\label{eqn: 2mL}
        \int_{T}^{2T}\Big|L\Big(\frac{1}{2}+it, \pi\Big)\Big |^2 \dd t\ll_{\pi} \frac{T^{\frac{d}{2}}}{\log^{\eta_d}T}
\end{equation}
for any $0\leq\eta_d\leq \frac{1}{400d^4}$.
\end{theorem}

Our theorem is closely related to the weak subconvexity bound for the automorphic $L$-functions of higher degree. There are many breakthrough regarding this topic. In \cite{Sound10}, assuming GRC, Soundararajan obtained a subconvex bound with a saving of one logarithm power relative to the analytic conductor of $L$-function. Later, in \cite{S-T19}, Soundararajan and Thorner achieved a small logarithmic saving without relying on GRC. Recently,  Nelson \cite{Nelson21} made a significant advance by proving a power-saving type subconvex bound for the standard $L$-function on the $t$-aspect.

The traditional approach to estimating the second moment of $L$-functions relies crucially on obtaining non-trivial bounds for the shifted convolution problem:
\[
\sum_{|h| \approx T^{d/2-1}}\sum_{n \approx T^{d/2}}\lambda_\pi(n)\overline{\lambda_\pi(n+h)}.
\]  
This sum captures the correlation of coefficients $\lambda_\pi(n)$ at scale $T^{d/2}$ at shifts of size about $T^{d/2-1}$, and non-trivial cancellation in such sums is essential for achieving power-saving bounds in moment estimates. When $\pi$ is of degree $3$, power-saving bounds for this shifted convolution can be achieved using the delta method, as demonstrated in \cite{Pal22}, \cite{DLY24}. 
However,  for automorphic representations $\pi$ of degree $d \geq 4$, the situation becomes significantly more difficult. The structure of the coefficients $\lambda_\pi(n)$ grows more complex, and the analytic tools, such as variations of the circle method or spectral methods, do not readily extend to produce non-trivial estimates. As a result, effective methods for controlling these shifted sums in higher degrees remain elusive, presenting a major obstacle to advancing the theory of high-degree $L$-functions.

For the proof of Theorem \ref{thm: 2M_L}, we rely on the multiplicative property of the Dirichlet coefficients $\lambda_{\pi}(n)$ and assume the Generalized Ramanujan Conjecture (GRC) in order to apply techniques from the theory of multiplicative functions.
Let $d\geq 3$ and $X\geq T^{\frac{d}{2}}(\log T)^{-A}$.
A central step in the argument is to establish the second moment estimate:
\begin{equation}\label{eqn: 2mlambda}
    \int_{T}^{2T}\Big|\sum_{n\leq X}\lambda_{\pi}(n)n^{-it}\Big|^2 \dd t\ll \frac{X^2}{\log^{\delta}X}
\end{equation}
for some $\delta>0$ (see Proposition \ref{prop: 2M_estimate}).
This is similar to the situation of Liouville function $\lambda(n)$ mentioned in \cite{M-R15, M-R16}. Indeed, Matom\"aki and Radziwi\l\l\ showed that
\begin{equation*}
  \int_{0}^{T}\Big|\sum_{X< n\leq 2X}\frac{\lambda(n)}{n^{1+it}}\Big|^2\dd t\ll \frac{1}{(\log X)^{\frac{1}{3}-o(1)}}\Big(\frac{T}{X}+1\Big)+\frac{T}{X^{1-o(1)}},
\end{equation*}
which is non-trivial for $X\geq T^{1+\varepsilon}.$
Following the ideas in \cite{M-R15} and their extension to the general divisor bounded multiplicative functions, when  $X\geq T^{1+\varepsilon},$ we find that
\begin{equation}\label{eqn:PretentiousLS}
 \int_{[T, 2T]\setminus B_{t_0}}\Big|\sum_{n\leq X}\lambda_{\pi}(n)n^{-it}\Big|^2 \dd t\ll \frac{X^2}{\log^{\delta}X},
\end{equation}
where $B_{t_0}\subset [T, 2T]$ is a possible small interval (may be of size $\log X$) such that $\lambda_{\pi}(n)$ may exhibit pretentious behavior toward $n^{it}$ with $t\in B_{t_0}$. On the other hand, under GRC, Soundararajan \cite{Sound10} established the pointwise upper bound estimate:
\begin{equation}\label{eqn:SoundBound}
 \sum_{n\leq X}\lambda_{\pi}(n)n^{-it}\ll_{\pi} \frac{X}{(\log X)^{1-\varepsilon}},
\end{equation}
uniformly for $X\gg t^{\frac{d}{2}}(\log t)^{-A}.$ 
Combining the above two bounds, we deduce the desired bound \eqref{eqn: 2mlambda} in the range of $X\geq T^{\frac{d}{2}}(\log T)^{-A}.$  
Both the method of Matom\"aki--Radziwi\l\l \ in proving \eqref{eqn:PretentiousLS} and Soundararajan's estimate \eqref{eqn:SoundBound} rely crucially on the assumption that the coefficients $\lambda_\pi(n)$ are divisor bounded. For this reason, our theorem ultimately depends on the assumption of GRC.

The second input is a combination of  approximate functional equation for $L$-functions and the partial summation. This gives
\[
L(\frac{1}{2}+it, \pi) =\int_{1}^{\infty}\sum_{n\leq x}\lambda_{\pi}(n)n^{-it}g(x,t)\frac{\dd x}{x^{3/2}}+\cdots
\]
where $g(x, t)$ satisfies the bound $g(x, t)\ll (\frac{t^{\frac{d}{2}}}{x})^{B}$ for any $B\gg 1.$ 
Here, it is crucial that the sum over $n$ is a discrete sum, which remains separate from the weight function. This separation ensures that the second moment estimate can be applied directly without interference from the smoothing weights. 
Then by using Minkowski's inequality, we have
\[
\int_{T}^{2T}\Big|L\Big(\frac{1}{2}+it, \pi\Big)\Big |^2\dd t \ll \left(\int_{1}^{\infty}\Big(\int_{T}^{2T}\Big|\sum_{n\leq x}\lambda_{\pi}(n)n^{-it}g(x,t)\Big|^2\dd t\Big)^{\frac{1}{2}}\frac{\dd x}{x^{3/2}}
\right)^2+\cdots.
\]
Using the mean value theorem trivially,  we get the $\varepsilon$-free bound
\[
\int_{T}^{2T}|L\Big(\frac{1}{2}+it, \pi\Big) |^2\dd t\ll T^{\frac{d}{2}}
\]
which matches the convexity bound for the second moment of $L$-functions. Further more, applying \eqref{eqn: 2mL} in certain range instead of the mean value theorem, we can obtain the log saving result for the second moment estimate.

If we assume GRC plus a strong zero free region for $L(s, \pi)$ of degree $d\geq 3$ as 
\begin{equation}\label{eqn: ZFR}
L(\sigma+it, \pi)\neq 0 \quad \text{ for } \sigma>1-\frac{c_\pi}{\log^{\beta}(3+|t|)}
\end{equation}
with constants $c_\pi>0$ and $0\leq \beta<1 $ which depends on $\pi$, then we can prove that
\[
\int_{T}^{2T}\Big|L\Big(\frac{1}{2}+it, \pi\Big)\Big |^2\dd t \ll_{\pi, \varepsilon}\frac{T^{\frac{d}{2}}}{\log^{1-\beta-\varepsilon}T}
\]
for any small $\varepsilon>0$. Indeed, under assumption \eqref{eqn: ZFR},  we can use a variant of prime number theorem instead of the pretentious approach for multiplicative functions in the proof of \eqref{eqn: 2mlambda}. The same argument is applied to the twisted $L$-functions on the character family in \cite{H-L24}.
Such a strong zero free region \eqref{eqn: ZFR} is extremely important in the theory of $L$-function.
It is only known for us in the case of Riemann $\zeta$ function and is out of reach by the current techniques in other automorphic  $L$-functions.

\medskip
\textbf{Notation.}
Throughout the paper, $\varepsilon$ is an arbitrarily small positive number and $A$ is an arbitrarily large positive number,
all of them may be different at each occurrence.
As usual, $p$ stands for a prime number.
We use the standard Landau and Vinogradov asymptotic notations \( O(\cdot) \), \( o(\cdot) \), \(\ll\), and \(\gg\). Specifically, we express \( X \ll Y \), \( X = O(Y) \), or \( Y \gg X \) when there exists a constant \( C \) such that \( |X| \leq CY \). If the constant $C=C_s$ depends on some object $s$, we write $X=O_s(Y)$. And the notation \( X \asymp Y \) is used when both \( X \ll Y \) and \( Y \ll X \) hold. As \( N \to \infty \), \( X = o(Y) \) indicates that \( |X| \leq c(N)Y \) for some function \( c(N) \) that tends to zero.
For a proposition \( P \), we write \( 1_P(x) \) for its indicator function, i.e. a function that equals to 1  if \( P \) is true and 0 if \( P \) is false.

\section{Background of $L$-functions}\label{sec:Background_L}
Let $d \geq 3$ and $\mathbb{A}$ be the ring of adel\'es of $\mathbb{Q}$, let $\pi$ be an irreducible cuspidal automorphic representation of $GL_d(\mathbb{A})$ with unitary central character that acts trivially on the diagonally embedded copy of $\mathbb{R}^{+}$. Let $q_\pi$ denote the conductor of $\pi$. The finite part of $\pi$ factors as a tensor product $\pi=\otimes_p \pi_p$, with local representations $\pi_p$ at each prime $p$. The local $L$-function at $p$ takes the form
$$
L\left(s, \pi_p\right)=\prod_{1 \leq j \leq d}\left(1-\frac{\alpha_{j, \pi}(p)}{p^s}\right)^{-1}:=\sum_{l \geq 0} \frac{\lambda_\pi\left(p^l\right)}{p^{l s}},
$$
where $\left\{\alpha_{1, \pi}(p), \cdots, \alpha_{d, \pi}(p)\right\} \subset \mathbb{C}$ are the Satake parameters of $\pi_p$. The standard $L$-function of $\pi$ is
$$
L(s, \pi):=\prod_p L\left(s, \pi_p\right)=\sum_{n \geq 1} \frac{\lambda_\pi(n)}{n^s},
$$
which converges absolutely when $\Re(s)>1$. The Dirichlet coefficient $\lambda_\pi(n)$  is multiplicative, with
\begin{equation}
    \lambda_\pi\left(p^r\right)=\sum_{\substack{r_1, \ldots, r_d \geq 0 \\ r_1+\cdots+r_d=r}} \prod_{1 \leq j \leq d} \alpha_{j, \pi}(p)^{r_j}.
\end{equation}
Especially,
$\lambda_\pi(p)=\sum_{1\leq j\leq d}\alpha_{j, \pi}(p).$
The generalized Ramanujan conjecture (GRC) asserts that for all $1 \leq j \leq d$,
\[
|\alpha_{j, \pi}(p)| \leq 1.
\]
It follows that if $\pi$ satisfies GRC then
$$
\left|\lambda_\pi\left(p^r\right)\right| \leq \sum_{\substack{r_1, \ldots, r_d \geq 0 \\
r_1+\cdots+r_d=r}} 1=\tau_d\left(p^r\right).
$$
Therefore $\left|\lambda_\pi(n)\right| \leq \tau_d(n)$.\par
Associated with $\pi$, there is also an Archimedean $L$-factors defined as
\begin{equation}
    L(s, \pi_{\infty}):={q_\pi}^{\frac{s}{2}}\pi^{-\frac{ds}{2}}\prod_{j=1}^{d}\Gamma\left(\frac{s+\mu_{j, \pi}}{2}\right),
\end{equation}
where $\mu_{j, \pi}$ are local parameters at $\infty$.
According to \cite{LRS}, \cite{KS03}, we have
\begin{equation}\label{eqn: theta_d}
    |\Re{\mu_{j, \pi}}|\leq \theta_d,
\end{equation}
where $\theta_2=\frac{7}{64} $,
$\theta_3=\frac{5}{14} $,
$\theta_4=\frac{9}{22} $ and
$\theta_d=\frac{1}{2}-\frac{1}{d^2+1}$ for $d\geq 5$.
This implies that $L(s, \pi_{\infty})$ is analytic in the right half plane $\Re(s)>\theta_d$.\par
The complete $L$-function
\begin{equation}
    \Lambda(s, \pi)=L(s, \pi)L(s, \pi_{\infty})
\end{equation}
admits the analytic continuation in the whole complex plane $\mathbb{C}$.
Moreover $\Lambda(s, \pi)$ satisfies a functional equation
\begin{equation}\label{eqn:func_eq}
    \Lambda(s, \pi)
    =\varepsilon_\pi
    \Lambda(s, \Tilde{\pi}),
\end{equation}
where $\varepsilon_\pi$ is a complex number of modulus $1$ and $\Tilde{\pi}$ is the contragredient of $\pi$. For any place $\ell\leq \infty$, $\Tilde{\pi}_\ell$ is equivalent to the complex conjugate $\overline{\pi}_\ell$, and we have
\begin{equation}
    \{\alpha_{j, \pi}(p)\}=\{\overline{\alpha}_{j, \pi}(p)\}, \quad
    \{\mu_{j, \pi}\}=\{\overline{\mu}_{j, \pi}\}.
\end{equation}
This implies $\lambda_{\pi}(n)=\overline{\lambda}_{\Tilde{\pi}}(n)$.\par
Taking the log-derivative of $L$-function, we can write, for $\Re(s)>1$,
\begin{equation}
    -\frac{L'}{L}(s,\pi)=\sum_{n\geq 1}\frac{\Lambda(n)a_{\pi}(n)}{n^s},
\end{equation}
where
\begin{equation}
    a_{\pi}(n)=
    \begin{cases}
        \sum_{1\leq j\leq d}\alpha_{j, \pi}(p)^k, &\text{if}\quad n=p^k,\\
         0, &\text{otherwise},\\
    \end{cases}
\end{equation}
and $\Lambda(n)$ is the von Mangoldt function.
Thus we have
\begin{equation}\label{eqn: lambda_p^k}
    \sum_{k\geq 0}\frac{\lambda_{\pi}(p^k)}{p^{ks}}=\exp\left(\sum_{v\geq 1}\frac{a_{\pi}(p^v)}{v}p^{-vs}\right).
\end{equation}
Especially $a_{\pi}(p)=\lambda_{\pi}(p)$ for any prime $p$.\par
We also introduce the Rankin--Selberg $L$-function $L(s, \pi \times \Tilde{\pi})$ associated to $\pi$ and $\Tilde{\pi}$. For unramified prime $p$, the local $L$-function takes the form
$$
L(s, \pi_p \times \Tilde{\pi}_p)=\prod_{1\leq j\leq d}\prod_{1\leq \ell\leq d}\left(1-\frac{\alpha_{j, \pi}\overline{\alpha}_{\ell, \Tilde{\pi}}}{p^s}\right)^{-1}:=\sum_{k=0}\frac{\lambda_{\pi\times \Tilde{\pi}}(p^k)}{p^{ks}},
$$
and the Rankin--Selberg $L$-function $L(s, \pi \times \Tilde{\pi})$ is defined by
$$
L(s, \pi \times \Tilde{\pi}):=\prod_p L\left(s, \pi_p \times \Tilde{\pi}_p\right)
=\sum_{n \geq 1} \frac{\lambda_{\pi\times \Tilde{\pi}}(n)}{n^s}
$$
which converges absolutely for $\Re(s)>1$.
Thus we have
\begin{equation}\label{eqn: lambda2_p^k}
    \sum_{k\geq 0}\frac{\lambda_{\pi \times \Tilde{\pi}}(p^k)}{p^{ks}}=\exp\left(\sum_{v\geq 1}\frac{|a_{\pi}(p^v)|^2}{v}p^{-vs}\right).
\end{equation}
From \eqref{eqn: lambda_p^k}, \eqref{eqn: lambda2_p^k} and \cite[Lemma 3.1]{Sound10}, we have, for an unramified prime $p$,
\begin{equation}
    |\lambda_{\pi}(p^k)|^2\leq \lambda_{\pi \times \Tilde{\pi}}(p^k),
\end{equation}
and
\begin{equation}\label{eqn: square_ineq}
    |\lambda_{\pi}(n)|^2\leq \lambda_{\pi \times \Tilde{\pi}}(n).
\end{equation}
In fact, by \cite[Lemma 3.1]{J-L-W21}, we know that \eqref{eqn: square_ineq} holds for all positive integer $n$.


We introduce the following various lemmas which are related to the coefficient $\lambda_{\pi}(n)$.
These lemmas are standard and will be used several times in our proof.

\begin{lemma}\label{lemma: weak_subconvexity}
    Let $|t|\geq 2$ and $|t|^{\frac{d}{2}}(\log |t|)^{-A}\leq X$, assuming GRC, we have
    \begin{equation}
        \sum_{n\leq X}\lambda_{\pi}(n)n^{-it}\ll_{\pi} \frac{X}{(\log X)^{1-\varepsilon}}.
    \end{equation}
\end{lemma}
\begin{proof}
    This is \cite[Theorem 2]{Sound10} on the $t$-aspect.
\end{proof}

The above estimate leads directly to a weak subconvexity bound for $L(1/2+it, \pi)$. The convexity bound for $L$-function give the same upper bound when $X\gg |t|^{\frac{d}{2}}(\log |t|)^{O(1)}$. Going beyond this range, specifically, achieving the same cancellation for values as small as 
$X=|t|^{\frac{d}{2}}(\log |t|)^{-A}$ is essential for obtaining subconvexity. Our second moment estimate \eqref{eqn: 2mlambda} is analogous to this situation.

\begin{lemma}\label{lemma: p_2moment}
    Let $X\geq 2$,  there is a constant $c>0$ depending on degree $d$ such that
    \begin{equation}
        \sum_{p\leq X}|\lambda_\pi(p)|^2\log p= X+O_{\pi}(Xe^{-c\sqrt{\log X}}).
    \end{equation}
    In particular, we have
    \begin{equation}
        \sum_{z< p \leq w}\frac{|\lambda_{\pi}(p)|^2}{p}=\sum_{z<p\leq w}\frac{1}{p}+O_{\pi}\big(\frac{1}{\log z}\big)
    \end{equation}
    for any $2\leq z\leq w$
\end{lemma}
\begin{proof}
    The proof is standard, see for example \cite[Lemma 5.6]{Alex23}.
\end{proof}

\begin{lemma}\label{lemma: coef_2moment}
    Let $X\geq 2$,  assuming GRC,  we have
    \begin{equation}
        \sum_{n\leq X}\lambda_{\pi\times\Tilde{\pi}}(n)=c_{\pi}X+O_{\pi,\varepsilon}(X^{\frac{d^2-1}{d^2+1}+\varepsilon}),
    \end{equation}
    for some constant $c_\pi$ depending on $\pi$.
    Note that $|\lambda_{\pi}(n)|^2\leq \lambda_{\pi\times\Tilde{\pi}}(n).$
    Therefore, for $X\geq Y\geq X^{\frac{d^2-1}{d^2+1}+\varepsilon}$, we have
    \begin{equation}
        \sum_{X\leq n\leq X+Y}|\lambda_{\pi}(n)|^2\ll_\pi Y.
    \end{equation}
\end{lemma}
\begin{proof}
    See e.g. \cite[Eqn (3.10)]{Lu09}, \cite[Proposition 1.1]{FI05}.
    \end{proof}
\begin{lemma}\label{lemma: sieve_bound}
    Let $X\geq 2$, and $1\ll P\leq Q\leq X$,  assuming GRC, we have
    \begin{equation}
            \sum_{n\leq X\atop (n, \prod_{P\leq p\leq Q}p)= 1}\lambda_{\pi\times\Tilde{\pi}}(n)\ll_\pi X\prod_{P\leq p\leq Q}\Big(1-\frac{1}{p}\Big).
    \end{equation}
Using $|\lambda_{\pi}(n)|^2\leq \lambda_{\pi\times\Tilde{\pi}}(n)$, we have
    \begin{equation}
            \sum_{n\leq X\atop (n, \prod_{P\leq p\leq Q}p)= 1}|\lambda_{\pi}(n)|^2\ll_\pi X\prod_{P\leq p\leq Q}\Big(1-\frac{1}{p}\Big).
    \end{equation}
\end{lemma}
\begin{proof}
    See \cite[Lemma 5.3]{J-L21}.
\end{proof}

The three lemmas above address various mean value estimates for the coefficients $\lambda_\pi(n)$. In general, under the assumption of GRC, 
the coefficients $\lambda_\pi(n)$ behave similarly to a bounded sequence on average, despite actually being divisor-bounded. Hence
when applying standard tools such as the mean value theorem to handle the second moment, summations involving  $\lambda_\pi(n)$ will not introduce some extraneous logarithmic factors. It is therefore advantageous to retain the slight saving established in \eqref{eqn: 2mlambda} throughout the argument. 

\section{Second moment for Dirichlet polynomials}
In this section, we present the proof of \eqref{eqn: 2mlambda}. 
Our approach follows an argument originally introduced by Matom\"aki and Radziwi\l\l \ \cite{M-R16} for real valued $1$-bounded multiplicative functions, developed to study multiplicative functions in short intervals.
This method was later extended to complex-valued $1$-bounded multiplicative functions in \cite{M-R20}, and further generalized by Mangerel \cite{Alex23} and Sun \cite{Sun24} to the divisor bounded multiplicative functions.
In the present setting, we apply this technique specifically to the coefficients $\lambda_\pi(n)$ and establish the following second moment estimate:

\begin{proposition}\label{prop: 2M_estimate}
Let $T\geq 2$ and $X\geq T^{\frac{d}{2}}(\log T)^{-A}$ where $d\geq 3$. Assuming GRC, then we have
    \begin{equation}
        \int_{2}^{T}\Big|\sum_{n\leq X}\lambda_{\pi}(n)n^{-it}\Big|^2 \dd t\ll \frac{X^2}{\log^{\eta_d}X},
    \end{equation}
    for any $0\leq \eta_d\leq \frac{1}{400d^4}$.
\end{proposition}

Proposition \ref{prop: 2M_estimate} provides a non-trivial improvement over the mean value theorem (Lemma \ref{lemma: mvt}),  which only yields the bound $O(X^2).$
The proof divides into two cases, depending on whether the coefficients $\lambda_{\pi}(n)$ behave like $n^{-it}$ for some $t\in [2, T]$ over range of $n\leq X.$
In the non-pretentious case, which accounts for most  $t\in [2, T]$, we apply the method of \cite{M-R16} to obtain the improved bound $O(X^2/\log^\delta X)$.
In the pretentious case, which occurs for a small exceptional set of $t$,  the pointwise bound from Lemma \ref{lemma: weak_subconvexity} suffices to complete the estimate. 

We first introduce the necessary tools in \S \ref{sec:3.1}, \S \ref{sec:3.2} and \S \ref{sec:3.3}. The proof of Proposition \ref{prop: 2M_estimate} is then presented in  \S \ref{sec:3.4}.

\subsection{Mean value theorems}\label{sec:3.1}
We introduce the following various mean value theorems for our subsequent estimate.
\begin{lemma}\label{lemma: mvt}
Let $\{a_n\}$ be a sequence of complex number,
   we have
   \begin{equation}\label{eqn:mvt}
       \int_{0}^{T}\Big|\sum_{N\leq n\leq 2N}a_nn^{-it}\Big|^2\dd t\ll (T+N)\sum_{N\leq n\leq 2N}|a_n|^2.
   \end{equation}
\end{lemma}
\begin{proof}
    See \cite[Lemma 6]{M-R16} or \cite[Theorem 9.1]{I-K04}.
\end{proof}

For the rest of the paper, we say that $\mathcal{T}\subseteq \mathbb{R}$ is well-spaced if $|t_1-t_2|\geq 1$ for all distinct $t_1, t_2\in \mathcal{T}.$ From the above mean value theorem, one can obtain the following large value result for Dirichlet polynomials.

\begin{lemma}\label{lemma: Tl_size}
   Let $\{a_p\}$ be a sequence of complex number with index prime $p$, satisfying $a_p\ll 1$,  and $\mathcal{T}\subset [-T, T]$ be a sequence of well-spaced points such that
   for every $t\in \mathcal{T}$,
   \[
   \Big|\sum_{P\leq p\leq 2P}a_pp^{it}\Big|\geq PV^{-1}
   \]
   where $V$ is a positive constant.
   Then
   \[
   |\mathcal{T}|\ll T^{2\frac{\log V}{\log P}}V^2\exp{\Big(2\frac{\log T}{\log P}\log\log T\Big)}
   \]
\end{lemma}
\begin{proof}
    See \cite[Lemma 8]{M-R16}.
\end{proof}

The following lemma is called the Hal\'asz inequality for primes in \cite{M-R16}.
\begin{lemma}\label{lemma: meanvalue_p}
    Let $\{a_p\}$ be a sequence of complex number with index prime $p$, and $\mathcal{T}\subset [-T, T]$ be a sequence of well-spaced points, then
    \begin{equation}
        \sum_{t\in \mathcal{T}}\big|\sum_{P\leq p \leq 2P}a_p p^{it}\big|^2\ll \left(P+|\mathcal{T}|\exp{\big(-\frac{\log P}{(\log T)^{2/3+\varepsilon}}\big)}(\log T)^2\right)\sum_{P\leq p\leq 2P}\frac{|a_p|^2}{\log P}.
    \end{equation}
\end{lemma}
\begin{proof}
    See \cite[Lemma 11]{M-R16}. The proof uses the duality principle and Korobov--Vinogradov type zero free region for the Riemann $\zeta$-function.
\end{proof}
\subsection{Decomposition of the Dirichlet polynomial}\label{sec:3.2}
The following lemma is a variant of \cite[Lemma 12]{M-R16}, which is used to construct the bilinear form for the second moment of Dirichlet polynomials.
 
\begin{lemma}\label{lemma: decop_2M}
    Let $H, T, X\geq 1$ and $Q\geq P\geq 1$. Let $f(n)$ be a $\tau_d$-bound multiplicative function and let $\mathcal{T}\subset [-T,T]$. Then, there exists some $j\in [\lfloor H\log P\rfloor, H\log Q]$ such that
    \begin{multline}\label{eqn:decop_2M}
        \int_{\mathcal{T}}\Big|\sum_{n\leq X}f(n)n^{-it}\Big|^2\dd t\ll \Big(H\log \frac{Q}{P}\Big)^2\int_{\mathcal{T}}\Big|Q_{j, H}(it)F_{j, H}(it)\Big|^2\dd t+
        (T+X)\frac{X\log^{d^2-1} X}{P}\\+
        (T+X)\sum_{P\leq p\leq Q}\sum_{X/p<m\leq Xe^{1/H}/p}|f(m)|^2+(T+X)\sum_{n\leq X\atop (n, \prod_{P\leq p\leq Q}p)= 1}|{f(n)}|^2.
    \end{multline}
    where
    \begin{equation}
        Q_{v, H}(s)=\sum_{P\leq p\leq Q\atop e^{v/H}\leq p \leq e^{(v+1)/H}}\frac{f(p)}{p^{s}},
    \end{equation}
    and
    \begin{equation}
        F_{v, H}(s)=\sum_{ m\leq Xe^{-v/H}}\frac{f(m)}{m^{s}(\omega_{[P,Q]}(m)+1)}.
    \end{equation}
\end{lemma}
\begin{proof}
    Writing $s=it$,  by Ramar\'e's identity
    \begin{equation}\label{eqn:Ramare}
        \sum_{mp=n\atop P\leq p\leq Q}\frac{1}{\omega_{[P, Q]}(m)+1_{p\nmid m}}=
        \begin{cases}
           1, &\quad \text{if} \quad (n, \prod_{P\leq p\leq Q}p)\neq 1,\\
           0, &\quad  \text{otherwise},\\
        \end{cases}
    \end{equation}
    where  $\omega_{[P, Q]}(m)$  is the number of distinct prime divisors of  $m$  in  $[P, Q]$,
    we have
    \begin{equation}\label{eqn: Dirichlet_f}
           \sum_{n\leq X}\frac{f(n)}{n^s}=\sum_{P\leq p\leq Q}\sum_{m\leq X/p}\frac{f(pm)}{(pm)^s(\omega_{[P, Q]}(m)+1_{p\nmid m})}+\sum_{n\leq X\atop (n, \prod_{P\leq p\leq Q}p)= 1}\frac{f(n)}{n^s}
    \end{equation}
    By the multiplicity of $f$, we rewrite the first sum as
    \begin{multline}
        \sum_{P\leq p\leq Q}\frac{f(p)}{p^s}\sum_{m\leq X/p}\frac{f(m)}{m^s(\omega_{[P, Q]}(m)+1)}
        \\+\sum_{P\leq p\leq Q}\sum_{m\leq X/p\atop p\mid m}\Big(\frac{f(pm)}{(pm)^s\omega_{[P, Q]}(m)}-\frac{f(p)f(m)}{(pm)^s(\omega_{[P, Q]}(m)+1)}\Big).
    \end{multline}
    We split the first sum above into dyadic ranges getting that
    \begin{equation}\label{eqn: decompostion_f}
        \sum_{\lfloor H\log P\rfloor\leq j\leq H\log Q}\sum_{e^{j/H}\leq p\leq e^{(j+1)/H)}\atop P\leq p\leq Q}\frac{f(p)}{p^s}\sum_{m\leq Xe^{-j/H}\atop mp\leq X}\frac{f(m)}{m^s(\omega_{[P,Q]}(m)+1)}.
    \end{equation}
    We drop the condition $mp\leq X$, which leads to a small overcount due to integers $mp$
    in the range $(X, Xe^{1/H}]$. Denote the  over-counted terms of $n=pm$ by $d_nn^{-s}$, we see that
     \begin{equation}
              |d_n|\leq\sum_{mp=n\atop  P\leq p\leq Q}\Big|\frac{f(p)f(m)}{\omega_{[P,Q]}(m)+1}\Big|.
     \end{equation}
     \eqref{eqn: decompostion_f} is equal to
      \begin{equation}
        \sum_{\lfloor H\log P\rfloor\leq j\leq H\log Q}\sum_{e^{j/H}\leq p\leq e^{(j+1)/H)}\atop P\leq p\leq Q}\frac{f(p)}{p^s}\sum_{m\leq Xe^{-j/H}}\frac{f(m)}{m^s(\omega_{[P,Q]}(m)+1)}-\sum_{X< n\leq Xe^{1/H}}\frac{d_n}{n^s}.
    \end{equation}
    Thus we can rewrite \eqref{eqn: Dirichlet_f} as
    \begin{multline}
          \sum_{\lfloor H\log P\rfloor\leq j\leq H\log Q}Q_{j, H}(s)F_{j, H}(s)
          +\sum_{P\leq p\leq Q}\sum_{m\leq X/p\atop p\mid m}\Big(\frac{f(pm)}{(pm)^s\omega_{[P, Q]}(m)}-\frac{f(p)f(m)}{(pm)^s(\omega_{[P, Q]}(m)+1)}\Big)\\-\sum_{X< n\leq Xe^{1/H}}\frac{d_n}{n^s}+\sum_{n\leq X\atop (n, \prod_{P\leq p\leq Q}p)= 1}\frac{f(n)}{n^s}
    \end{multline}
    We square the above formula,  integrate over $\mathcal{T}$, and take the maximum term of the first sum over $j$ and apply the mean value theorem (Lemma \ref{lemma: mvt}) to the remaining term. The second sum contributes at most
    \begin{equation}
        (T+X)\sum_{mp^2\leq X\atop P\leq p\leq Q}\tau_d^2(mp^2)\ll (T+X)\sum_{P\leq p\leq Q}\sum_{m\leq X/p^2}\tau_d^2(m)\ll (T+X)\frac{X\log^{d^2-1} X}{P}.
    \end{equation}
    The third sum contributes at most
       \begin{equation*}
       \begin{split}
        (T+X)\sum_{X< n\leq Xe^{1/H}}&\Big(\sum_{mp=n\atop  P\leq p\leq Q}\Big|\frac{f(p)f(m)}{\omega_{[P,Q]}(m)+1}\Big|\Big)^2\\
        &\ll (T+X)\sum_{X< n\leq Xe^{1/H}}\omega_{[P, Q]}(n)^2\sum_{mp=n\atop  P\leq p\leq Q}\Big|\frac{f(p)f(m)}{\omega_{[P,Q]}(m)+1}\Big|^2\\
        &\ll(T+X)\sum_{X< n\leq Xe^{1/H}}\sum_{n=pm\atop P\leq p\leq Q}|f(m)|^2\\
        &\ll (T+X)\sum_{P\leq p\leq Q}\sum_{X/p<m\leq Xe^{1/H}/p}|f(m)|^2.
       \end{split}
      \end{equation*}
     The last sum contributes at most
       \begin{equation}
           (T+X)\sum_{n\leq X\atop (n, \prod_{P\leq p\leq Q}p)= 1}|{f(n)}|^2.
       \end{equation}
     Lemma \ref{lemma: decop_2M} is followed upon combining the above bounds.
\end{proof}
We will take $f(n)=\lambda_\pi(n)$ in Lemma \ref{lemma: decop_2M}.  
To bound the last two terms of \eqref{eqn:decop_2M} in upper bounds,
we use the mean value lemmas in \S \ref{sec:Background_L} instead of divisor bounds for $f$.

\subsection{Hal\'asz type estimates}\label{sec:3.3}
In this subsection, we deduce some Hal\'asz type estimates for $\lambda_{\pi}(n)$.  From now on, we always assume that GRC holds for $\pi$. Hence $\lambda_{\pi}$ is a $\tau_d$-bounded multiplicative function.
We introduce some tools about the pretentious approach in \cite{M-R20}. 
To deal with the multiplicative functions of complex value, we need the following definition.
\begin{definition}
Let $f:\mathbb{N}\rightarrow \mathbb{C}$ be a multiplicative function. We define $t_0:=t_0(f, X)\in [-X, X]$ to be (one of) $t\in [-X, X]$ which minimizes the following distance
\[
\mathbb{D}(f, n^{it}; X)^2:=\sum_{p\leq X}\frac{|f(p)|-\Re(f(p)p^{-it})}{p}.
\]
And we denote the minimum by
\[
 M(f; X)=\min_{|t|\leq X}\sum_{p\leq X}\frac{|f(p)|-\Re{(\lambda_{\pi}(p)p^{-it})}}{p}
\]
\end{definition}
Throughout the paper, 
we set $t_0\in [-X,X]$ (depending on $\pi$ and $X$) such that $\mathbb{D}(\lambda_\pi, n^{t_0}; X)^2$
attains the minimum $M(\lambda_\pi; X).$

\begin{lemma}\label{lemma: HalaM}
    For all $t\in [-X, X]$,
    \begin{equation}
         M(\lambda_{\pi}; X)\geq \left(\frac{1}{3d}\big(1-\frac{2d^2}{\pi}\sin{\frac{\pi}{2d^2}}\big)-\varepsilon\right)\min\{\log\log X, \log(1+|t-t_0|\log X)\}+O_{\pi}(1),
    \end{equation}
  for any small $\varepsilon>0.$
\end{lemma}
\begin{proof}
    By Lemma \ref{lemma: p_2moment} and GRC, using Cauchy--Schwarz, we have, for all $2\leq z\leq w\leq X,$
    \begin{equation}\label{eqn: lowbound_p}
        \sum_{z< p \leq w}\frac{|\lambda_{\pi}(p)|}{p}\geq \frac{1}{d}\sum_{z<p\leq w}\frac{1}{p}+O_{\pi}\big(\frac{1}{\log z}\big).
    \end{equation}
    Therefore,  this lemma follows from a consequence of \cite[Lemma 5.1]{M-R20} with $f(p)=\frac{\lambda_{\pi}(p)}{d}$.
\end{proof}

From now on, we define $\rho_d>0$ is a constant which is strictly less than $\frac{1}{3d}\big(1-\frac{2d^2}{\pi}\sin{\frac{\pi}{2d^2}}\big)$, for example $\rho_d=\frac{1}{100d^3}$.
Combining the above lemma with  Hal\'asz's theorem, we get the following estimate.

\begin{lemma}\label{lemma: Hala_bound1}
    Let $X^{1/5}\leq Y\leq X$, and $2\leq P\leq Q\leq \exp(\frac{\log X}{\log\log X})$. Then for any $1\leq Z\leq \log X$, we have
    \begin{equation*}
        \sup_{Z\leq |u|\leq X/2}\Big|\sum_{n\leq Y\atop p\mid n\implies p\notin [P, Q]}\lambda_\pi(n)n^{-i(t_0+u)}\Big|\ll Y\left(\Big(\frac{\log Q}{\log P}\Big)^{2d}\frac{\log\log X}{(\log X)^{\rho_d}}+\frac{1}{\sqrt{Z}}\right).
    \end{equation*}
\end{lemma}
\begin{proof}
    Under the assumption of GRC, using \eqref{eqn: lowbound_p} and Lemma \ref{lemma: HalaM}, this lemma is a consequence of \cite[Lemma 3.3]{Alex23} with $A=\frac{1}{d}, B=d, C=1, \gamma=1$ and $\sigma=\rho_d$.  Note that by Lemma \ref{lemma: p_2moment}, using Cauchy--Schwarz, we have, for $X\gg 1$,
    \[
    \sum_{p\leq X}\frac{|\lambda_\pi(p)|}{p}\ll \Big(\sum_{p\leq X}\frac{|\lambda_\pi(p)|^2}{p}\Big)^{\frac{1}{2}}\Big(\sum_{p\leq X}\frac{1}{p}\Big)^{\frac{1}{2}}\ll \log\log X\sim \sum_{p\leq X}\frac{1}{p}.
    \]
    This implies that, in \cite[Lemma 3.3]{Alex23}, 
    $$\mathcal{P}_{\lambda_\pi}(X):=\prod_{p\leq X}\Big(1+\frac{|\lambda_\pi(p)|-1}{p}\Big)\ll \exp\Big(\sum_{p\leq X}\frac{|\lambda_\pi(p)|-1}{p}\Big)\ll 1$$
    which gives the desired upper bound.
\end{proof}

In the proof of Proposition \ref{prop: 2M_estimate}, we need to apply a Hal\'asz's type estimate to a function which is not quite multiplicative and the following lemma takes care of this application to a polynomial arising from the Ramar\'e's identity \eqref{eqn:Ramare}.

\begin{lemma}\label{lemma: Hala_bound2}
    Let $X>100$ large, $X^{1/2}\leq Y\leq X$, $2\leq P\leq Q\leq \exp(\frac{\log X}{\log\log X})$ and $1\leq Z\leq \log X$, we have
    \begin{equation*}
        \sup_{Z\leq u\leq X/2}\Big|\sum_{ n\leq Y}\frac{\lambda_\pi(n)}{n^{i(t_0+u)}(\omega_{[P, Q]}(n)+1))}\Big|\ll Y\left(\Big(\frac{\log Q}{\log P}\Big)^{3d}\frac{\log\log X}{(\log X)^{\rho_d}}+\Big(\frac{\log Q}{\log P}\Big)^{d}\frac{1}{\sqrt{Z}}\right).
    \end{equation*}
\end{lemma}
\begin{proof}
   Writing $n=n_1n_2$ where $n_1$ has all prime factors from $[P, Q]$ and $n_2$ has none, we have
 \begin{multline*}
       \sum_{n\leq Y}\frac{\lambda_\pi(n)}{n^{i(t_0+u)}(\omega_{[P, Q]}(n)+1))}=\sum_{n_1\leq Y^{1/2}\atop p\mid n_1\implies p\in [P, Q]}\frac{\lambda_{\pi}(n_1)}{{n_1}^{i(t_0+u)}(\omega_{[P, Q]}(n_1)+1))}\sum_{n_2\leq Y/n_1\atop p\mid n_2\implies p\notin [P, Q]}
       \frac{\lambda_\pi(n_2)}{{n_2}^{i(t_0+u)}}\\
       +O\left(\sum_{n_2\leq Y^{1/2}\atop p\mid n_2\implies p\notin [P, Q]}|\lambda_\pi(n_2)|\sum_{n_1\leq Y/n_2\atop p\mid n_1\implies p\in [P, Q]}|\lambda_\pi(n_1)|\right)\\
       \ll \sum_{n_1\leq Y^{1/2}\atop p\mid n_1\implies p\in [P, Q]}\tau_d(n_1)\Big|\sum_{n_2\leq Y/n_1\atop p\mid n_2\implies p\notin [P, Q]}
       \frac{\lambda_\pi(n_2)}{{n_2}^{i(t_0+u)}}\Big|
       +\sum_{n_2\leq Y^{1/2}\atop p\mid n_2\implies p\notin [P, Q]}\tau_d(n_2)\sum_{n_1\leq Y/n_2\atop p\mid n_1\implies p\in [P, Q]}\tau_d(n_1).
   \end{multline*}
For the first term, we use Lemma \ref{lemma: Hala_bound1} ($Y/n_1\geq Y^{1/2}/2\geq X^{1/5}$) for the inner sum and Shiu's theorem \cite[Theorem 1]{Shiu1980} for the outer sum. It is bounded by
$$Y\left(\Big(\frac{\log Q}{\log P}\Big)^{3d}\frac{\log\log X}{(\log X)^{\rho_d}}+\Big(\frac{\log Q}{\log P}\Big)^{d}\frac{1}{\sqrt{Z}}\right).$$
For the second term, using the $Q$-smooth estimate for $\tau_d$ (for example, see \cite[Lemma 3.4]{Alex23}),
we have,  the inner sum is bounded by
\begin{equation*}
    \sum_{ n_1\leq Y/n_2\atop P_{+}(n_1)\leq Q}\tau_d(n_1)\ll\frac{Y}{n_2} (\log Q)^{d-1}\exp\big(-\frac{1}{3} y\log y\big)\ll \frac{Y}{n_2(\log X)^{2A}},
\end{equation*}
where $y=\frac{\log (Y/n_2)}{\log Q}\geq \frac{1}{2}\log\log X$.
Thus the second term is bounded by
\begin{equation*}
    \frac{Y}{(\log X)^{2A}}\sum_{n_2\leq Y^{1/2}\atop p\mid n_2\implies p\notin [P, Q]}\frac{\tau_d(n_2)}{n_2}\ll\frac{Y}{(\log X)^{A}}.
\end{equation*}
Combining the above two bounds, we finish the proof.
\end{proof}

\subsection{Proof of Proposition \ref{prop: 2M_estimate}}\label{sec:3.4}
Let $B_Z(t_0):=[t_0-Z, t_0+Z]$ be an interval near the minimizer $t_0$ (depending on $\pi$ and $X$,
we split the integral range of $t$ int $B_Z(t_0)$ and $[2, T]\setminus B_Z(t_0).$
By Lemma \ref{lemma: weak_subconvexity}, we have
\begin{multline}\label{eqn: 2M_split}
    \int_{2}^{T}\Big|\sum_{n\leq X}\lambda_{\pi}(n)n^{-it}\Big|^2 \dd t \leq \int_{[2,T]\setminus B_{Z}(t_0)}\Big|\sum_{n\leq X}\lambda_{\pi}(n)n^{-it}\Big|^2 \dd t +2Z\sup_{t\in B_Z(t_0)}|\sum_{n\leq X}\lambda_{\pi}(n)n^{-it}\Big|^2\\
    \ll \int_{[2,T]\setminus B_{Z}(t_0)}\Big|\sum_{n\leq X}\lambda_{\pi}(n)n^{-it}\Big|^2 \dd t +\frac{ZX^2}{(\log X)^{2-\varepsilon}}.
\end{multline}
Using Lemma \ref{lemma: decop_2M}, we have, for  some $j\in [\lfloor H\log P\rfloor, H\log Q]$,
\begin{multline}\label{eqn: decomp_bound}
     \int_{[2,T]\setminus B_{Z}(t_0)}\Big|\sum_{n\leq X}\lambda_{\pi}(n)n^{-it}\Big|^2 \dd t \ll \big(H\log \frac{Q}{P}\Big)^2\int_{[0, T]\setminus B_{Z}(t_0)}\Big|Q_{j, H}(it)F_{j, H}(it)\Big|^2\dd t\\
     +(T+X)\frac{X\log^{d^2-1} X}{P}
        +(T+X)\sum_{P\leq p\leq Q}\sum_{X/p<m\leq Xe^{1/H}/p}|\lambda_{\pi}(m)|^2\\
        +(T+X)\sum_{n\leq X\atop (n, \prod_{P\leq p\leq Q}p)= 1}|{\lambda_{\pi}(n)}|^2.
\end{multline}
where $Q_{j, H}, F_{j, H}$ is defined in Lemma \ref{lemma: decop_2M} replacing $f$ by $\lambda_{\pi}.$
In the above formula, we set the parameters as follow:
\begin{equation}\label{eqn: parameters}
\begin{split}
    &Q=\exp{\big(\frac{\log X}{\log\log X}\big)},
    \quad\
    P=\exp{(\log^{1-\frac{\rho_d}{(3d+2)}}X)},\\
    &H=(\log X)^{\frac{\rho_d}{3d+2}},
    \quad\quad\quad
    Z=\log X.
\end{split}
\end{equation}
The second term in \eqref{eqn: decomp_bound} is bounded by  $O(X^2{\log^{-A} X})$ trivially.
For the third term in \eqref{eqn: decomp_bound}, using Lemma \ref{lemma: coef_2moment}, it is bound by
\[
    (T+X)\sum_{P\leq p\leq Q}\frac{X}{Hp}\ll (T+X)\frac{X}{H}\log\log Q\ll \frac{X^2\log\log X}{(\log X)^{\frac{\rho_d}{3d+2}}}.
\]
For the last term in \eqref{eqn: decomp_bound}, using Lemma \ref{lemma: sieve_bound}, it is bound by
\[
(TX+X^2)\frac{\log P}{\log Q}\ll\frac{X^2\log\log X}{(\log X)^{\frac{\rho_d}{3d+2}}}.
\]
Next, we give the estimate for the first term in \eqref{eqn: decomp_bound}.  Define
\[
\mathcal{T}_S:=\{t\in [2,T]\setminus B_{Z}(t_0): |Q_{j, H}(it)|\leq e^{(j+1)/H}(\log X)^{-100}\}
\]
and
\[
\mathcal{T}_L:=\{t\in [2,T]\setminus B_{Z}(t_0): |Q_{j, H}(it)|>e^{(j+1)/H}(\log X)^{-100}\}.
\]
On $\mathcal{T}_S$, by Lemma \ref{lemma: mvt}, we have
\begin{multline}\label{eqn: Ts_int}
    \big(H\log \frac{Q}{P}\Big)^2\int_{\mathcal{T}_S}\Big|Q_{j, H}(it)F_{j, H}(it)\Big|^2\dd t\\
    \ll H^2(\log Q)^2(\log X)^{-200}e^{(2j+2)/H}\int_{[0, T]}\Big|F_{j, H}(it)\Big|^2\dd t\\
    \ll H^2(\log Q)^2(\log X)^{-200}e^{(2j+2)/H}(T+Xe^{-j/H})Xe^{-j/H}\\
    \ll (TX^{1+\varepsilon}+X^2)(\log X)^{-100}\ll X^2(\log X)^{-100} .
\end{multline}
Let us now turn on $\mathcal{T}_L$.  We can find a well-spaced subset $\mathcal{T}_\ell\subset \mathcal{T}_L$ such that
\begin{equation}\label{eqn: Tl_int}
    \big(H\log \frac{Q}{P}\Big)^2\int_{\mathcal{T}_L}\Big|Q_{j, H}(it)F_{j, H}(it)\Big|^2\dd t\ll H^2(\log Q)^2\sum_{t\in \mathcal{T}_\ell}|Q_{j, H}(it)F_{j, H}(it)|^2
\end{equation}
By Lemma \ref{lemma: Tl_size}, we have
\begin{multline}\label{eqn: T_size}
    |\mathcal{T}_\ell|\ll \exp{\left(200\frac{\log\log X}{j/H}\log T+200\log\log X+2\frac{\log T\log\log T}{j/H}\right)}\\
    \ll  \exp{\left(\frac{(\log X)^{1+o(1)}}{\log P}\right)}\ll \exp{\left((\log X)^{\frac{\rho_d}{3d+2}+o(1)}\right)}.
\end{multline}
Using Lemma \ref{lemma: Hala_bound2}, we have
\begin{multline}\label{eqn: supF_it}
     \sup_{t\in [0,T]\setminus B_Z(t_0)}|F_{j, H}(it)|\ll Xe^{-j/H}\left(\Big(\frac{\log Q}{\log P}\Big)^{3d}\frac{\log\log X}{(\log X)^{\rho_d}}+\Big(\frac{\log Q}{\log P}\Big)^d\frac{1}{\sqrt{Z}}\right)\\
     \ll Xe^{-j/H}\Big(\frac{\log\log X}{(\log X)^{\frac{2\rho_d}{3d+2}}}+\frac{1}{(\log X)^{\frac{1}{2}-\frac{d\rho_d}{3d+2}}}\Big).
\end{multline}
Thus by Lemma \ref{lemma: meanvalue_p} with \eqref{eqn: T_size} and \eqref{eqn: supF_it}, we see that \eqref{eqn: Tl_int} is bounded by
\begin{equation*}
\begin{split}
     &H^2(\log Q)^2X^2e^{-2j/H}\Big(\frac{(\log\log X)^2}{(\log X)^{\frac{4\rho_d}{3d+2}}}+\frac{1}{(\log X)^{1-\frac{2d\rho_d}{3d+2}}}\Big)\sum_{t\in \mathcal{T}_\ell}|Q_{j, H}(it)|^2\\
     &\ll H^2(\log Q)^2X^2e^{-2j/H}\Big(\frac{(\log\log X)^2}{(\log X)^{\frac{4\rho_d}{3d+2}}}+\frac{1}{(\log X)^{1-\frac{2d\rho_d}{3d+2}}}\Big) \\
     &\times \left(e^{(j+1)/H}+\exp{\Big((\log X)^{\frac{\rho_d}{3d+2}+o(1)}-\frac{j+1}{H(\log T)^{\frac{2}{3}+\varepsilon}}}\Big)(\log T)^2\right)\sum_{e^{j/H}\leq p< e^{(j+1)/H}}|\lambda_{\pi}(p)|^2\frac{H}{j}.
\end{split}
\end{equation*}
Applying the Brun--Titchmarch inequality, we have
\begin{equation}
    \sum_{e^{j/H}\leq p< e^{(j+1)/H}}|\lambda_{\pi}(p)|^2\ll \sum_{e^{j/H}\leq p< e^{(j+1)/H}}1\ll \frac{e^{j/H}}{j}.
\end{equation}
Note that $\frac{j+1}{H(\log T)^{\frac{2}{3}+\varepsilon}}\geq \frac{\log P}{(\log X)^{\frac{2}{3}+\varepsilon}}\geq (\log X)^{\frac{\rho_d}{3d+2}+o(1)}$ and $j\geq H\log P$, we get that
\eqref{eqn: Tl_int} is bounded by
\begin{multline}
      H\big(\frac{\log Q}{\log P}\big)^2X^2\left(\frac{(\log\log X)^2}{(\log X)^{\frac{4\rho_d}{3d+2}}}+\frac{1}{(\log X)^{1-\frac{2d\rho_d}{3d+2}}}\right) \\
      \ll X^2\left(\frac{(\log\log X)^2}{(\log X)^{\frac{\rho_d}{3d+2}}}+\frac{1}{(\log X)^{1-\frac{(2d-3)\rho_d}{3d+2}}}\right)
      \ll X^2\frac{(\log\log X)^2}{(\log X)^{\frac{\rho_d}{3d+2}}}.
\end{multline}
Combining with \eqref{eqn: Ts_int}, we have
\begin{equation}
    \big(H\log \frac{Q}{P}\Big)^2\int_{[2, T]\setminus B_{Z}(t_0)}\Big|Q_{j, H}(it)F_{j, H}(it)\Big|^2\dd t
    \ll X^2\frac{(\log\log X)^2}{(\log X)^{\frac{\rho_d}{3d+2}}}.
\end{equation}
Collecting the above bounds, we have
\begin{equation}
    \int_{[2,T]\setminus B_{Z}(t_0)}\Big|\sum_{n\leq X}\lambda_{\pi}(n)n^{-it}\Big|^2 \dd t \ll
    X^2\frac{(\log\log X)^2}{(\log X)^{\frac{\rho_d}{3d+2}}}.
\end{equation}
Finally, the contribution for the integral on the small interval $B_Z(t_0)$ in \eqref{eqn: 2M_split} contributes at most
$
O(\frac{ZX^2}{(\log X)^{2-\varepsilon}})=O(\frac{X^2}{(\log X)^{1-\varepsilon}}).
$
This completes the proof of Proposition \ref{prop: 2M_estimate}.

\section{Proof of Theorem \ref{thm: 2M_L}}
Let $t\geq 2$, we shall adopt a standard argument about the approximate functional equation for $L(\frac{1}{2}+it, \pi)$,
see for example  Harcos \cite[Theorem 2.1]{Harcos} and Soundararajan \cite[Page 1486]{Sound10}.
For $c\geq \frac{1}{2},$  we move the line of integration
\[
\frac{1}{2\pi i}\int_{(c)}L(s+1/2+it, \pi)\frac{L(s+\frac{1}{2}+it, \pi_{\infty})}{L(\frac{1}{2}+it, \pi_{\infty})}e^{s^2}\frac{\dd s}{s}
\]
to the vertical line $\Re(s)=-c$, We encounter a pole at $s=0$,and so the above equals
\[
L(1/2+it, \pi)+\frac{1}{2\pi i}\int_{(-c)}L(s+1/2+it, \pi)\frac{L(s+\frac{1}{2}+it, \pi_{\infty})}{L(\frac{1}{2}+it, \pi_{\infty})}e^{s^2}\frac{\dd s}{s}.
\]
Now we use the functional equation \eqref{eqn:func_eq}, and make a change of variables $s\to -s$.
In this way we obtain that
\begin{multline}\label{eqn: FE}
      L(1/2+it, \pi)=\frac{1}{2\pi i}\int_{(c)}L(s+1/2+it, \pi)\frac{L(s+\frac{1}{2}+it, \pi_{\infty})}{L(\frac{1}{2}+it, \pi_{\infty})}e^{s^2}\frac{\dd s}{s}\\
      +\frac{\varepsilon_\pi}{2\pi i}\int_{(c)}L(s+1/2-it, \Tilde{\pi})\frac{L(s+\frac{1}{2}-it, \Tilde{\pi}_{\infty})}{L(\frac{1}{2}+it, \pi_{\infty})}e^{s^2}\frac{\dd s}{s}.
\end{multline}
Consider the first integral above, the second can be estimated similarly.
Using partial summation formula, we get
\begin{equation}
    L(s+1/2+it, \pi)=(s+1/2)\int_{1}^{\infty}\sum_{n\leq x}\lambda_{\pi}(n)n^{-it}\frac{\dd x}{x^{s+3/2}}.
\end{equation}
Insert it into the approximate functional equation \eqref{eqn: FE}, and exchanging the order of integration, we see that the first integral in \eqref{eqn: FE} equals
\begin{equation}\label{eqn: int_form}
    \int_{1}^{\infty}\sum_{n\leq x}\lambda_{\pi}(n)n^{-it}g(x,t)\frac{\dd x}{x^{3/2}},
\end{equation}
where
\begin{equation}
    g(x,t)=\frac{1}{2\pi i}\int_{(c)}(s+1/2)\frac{L(s+\frac{1}{2}+it, \pi_{\infty})}{L(\frac{1}{2}+it, \pi_{\infty})}e^{s^2}x^{-s}\frac{\dd s}{s}.
\end{equation}
 We can move the integral line to any $\Re(s)>\theta_d-\frac{1}{2}$, using Stirling's formula, we have
\begin{equation}\label{eqn: g-bound}
    g(x,t)\ll_{B,\pi}\Big(\frac{t^{\frac{d}{2}}}{x}\Big)^{B}, \quad \text{ for any } B\geq \theta_d-\frac{1}{2}+\varepsilon,
\end{equation}
where $\theta_d$ is introduced \eqref{eqn: theta_d}.

By the elementary inequality $(a+b)^2\leq 2(a^2+b^2)$, it suffices to bound the second moment of two integral in \eqref{eqn: FE}.
Here we only consider the second moment of the first integral since that the second can be estimated similarly.
By \eqref{eqn: int_form} and Minkowski's inequality, we have
\begin{equation*}
\begin{split}
    \int_{T}^{2T}\Big|\frac{1}{2\pi i}\int_{(c)}L(s&+1/2+it, \pi)\frac{L(s+\frac{1}{2}+it, \pi_{\infty})}{L(\frac{1}{2}+it, \pi_{\infty})}e^{s^2}\frac{\dd s}{s}\Big|^2\dd t\\
    &=\int_{T}^{2T}\Big|\int_{1}^{\infty}\sum_{n\leq x}\lambda_{\pi}(n)n^{-it}g(x,t)\frac{\dd x}{x^{3/2}}\Big|^2\dd t\\
    &\ll \left(\int_{1}^{\infty}\Big(\int_{T}^{2T}\Big|\sum_{n\leq x}\lambda_{\pi}(n)n^{-it}g(x,t)\Big|^2\dd t\Big)^{\frac{1}{2}}\frac{\dd x}{x^{3/2}}
    \right)^2\\
    &\ll \left(\int_{1}^{T^{d/2}(\log T)^{-A}}\Big(\int_{T}^{2T}\Big|\sum_{n\leq x}\lambda_{\pi}(n)n^{-it}g(x,t)\Big|^2\dd t\Big)^{\frac{1}{2}}\frac{\dd x}{x^{3/2}}\right)^2\\
    &\hspace{5em}+\left(\int_{T^{d/2}(\log T)^{-A}}^{\infty}\Big(\int_{T}^{2T}\Big|\sum_{n\leq x}\lambda_{\pi}(n)n^{-it}g(x,t)\Big|^2\dd t\Big)^{\frac{1}{2}}\frac{\dd x}{x^{3/2}}
    \right)^2.
\end{split}
\end{equation*}
The first part, using Lemma \ref{lemma: mvt}, Lemma \ref{lemma: coef_2moment} and \eqref{eqn: g-bound}, contributes at most
\begin{equation*}
\begin{split}
\Big(\int_{1}^{T^{d/2}(\log T)^{-A}}&\Big((T+x)x\Big)^{\frac{1}{2}}\Big(\frac{T^{\frac{d}{2}}}{x}\Big)^{B}\frac{\dd x}{x^{3/2}}\Big)^2\\
&\ll  T^{dB+1}\left(\int_{1}^{T^{d/2}(\log T)^{-A}}\frac{\dd x}{x^{B+1}}\right)^2+T^{dB}\left(\int_{1}^{T^{d/2}(\log T)^{-A}}\frac{\dd x}{x^{B+1/2}}\right)^2\\
&\ll T^{d/2}(\log T)^{-\frac{A}{2}},
\end{split}
\end{equation*}
by taking $B=1/4$.
The second part, we split the integral of $x$ into $[T^{d/2}(\log T)^{-A}, T^{d/2}) $ and $[T^{d/2}, +\infty).$
Using Proposition \ref{prop: 2M_estimate} and \eqref{eqn: g-bound}, it is bounded by
\begin{equation*}
\begin{split}
\Big(\int_{T^{d/2}(\log T)^{-A}}^{T^{d/2}}&\Big(\frac{x^2}{\log^{\eta_d}x}\Big)^{\frac{1}{2}}\Big(\frac{T^{\frac{d}{2}}}{x}\Big)^{B_1}\frac{\dd x}{x^{3/2}}\Big)^2
+\Big(\int_{T^{d/2}}^{\infty}\Big(\frac{x^2}{\log^{\eta_d}x}\Big)^{\frac{1}{2}}\Big(\frac{T^{\frac{d}{2}}}{x}\Big)^{B_2}\frac{\dd x}{x^{3/2}}\Big)^2
\\
&\ll
 \frac{T^{dB_1}}{\log^{\eta_d}T}
 \left(\int_{T^{d/2}(\log T)^{-A}}^{T^{d/2}}\frac{\dd x}{x^{B_1+\frac{1}{2}}}\right)^2
 + \frac{T^{dB_2}}{\log^{\eta_d}T}
 \left(\int_{T^{d/2}}^{{\infty}}\frac{\dd x}{x^{B_2+\frac{1}{2}}}\right)^2\\
&\ll \frac{T^{\frac{d}{2}}}
 {\log^{\eta_d}T}
\end{split}
\end{equation*}
by taking $B_1=1/3$ in first integral and $B_2=2$ in second integral. Theorem \ref{thm: 2M_L} follows upon combining the above bounds.

\section*{Acknowledgements}
The author wants to thank Prof. Bingrong Huang for his advices and encouragements.
He also wants to thank Haozhe Gou for many helpful discussions on \cite{M-R16}.
The author would like to thank the referees for their very helpful comments and suggestions.

\textbf{Conflict of interest.} The  author states that there is no conflict of interest.

\textbf{Data Availability.} Data sharing not applicable to this article as no datasets were generated or analysed during the current study.



\begin{thebibliography}{10}

\bibitem{CFKRS}
Conrey, John Brian; Farmer, David W.; Keating, Jonathan P.; Rubinstein, Michael O.; Snaith, Nina Claire. 
Integral moments of $L$-functions. 
\emph{Proc. London Math. Soc.} (3) 91 (2005), no. 1, 33--104.

\bibitem{DLY24}
Dasgupta, Agniva; Leung, Wing Hong; Young, Matthew Patrick.
The second moment of the $\rm GL_3$ standard $L$-function on the critical line.
\emph{ArXiv preprints}
arXiv:2407.06962

\bibitem{FI05}
Friedlander, John B.; Iwaniec, Henryk.
Summation formulae for coefficients of $L$-functions. 
\emph{Canad. J. Math.} 57 (2005), no. 3, 494--505. 


\bibitem{Go82}
Good, Anton.
The square mean of Dirichlet series associated with cusp forms.
\emph{Mathematika} 29 (1982), no. 2, 278--295 (1983).

\bibitem{H-L24}
Gou, Haozhe; Li, Liangxun.
On the second moment of twisted higher degree $L$-functions.
\emph{ArXiv preprints}
arXiv:2410.20144 

\bibitem{Harcos}
Harcos, Gergely.
Uniform approximate functional equation for principal $L$-functions.
\emph{Int. Math. Res. Not.} 2002, no. 18, 923--932.

\bibitem{Hardy-Littlewood}
Hardy, Godfrey Harold; Littlewood, John Edensor. 
Contributions to the theory of the riemann zeta-function and the theory of the distribution of primes. 
\emph{Acta Math.} 41 (1916), no. 1, 119--196.
 
\bibitem{Harper13}
Harper, Adam J.
Sharp conditional bounds for moments of the Riemann zeta function
\emph{ArXiv preprints} arXiv:1305.4618

\bibitem{HB}
Heath-Brown, D. R. 
The twelfth power moment of the Riemann-function. 
\emph{Quart. J. Math. Oxford Ser.} (2) 29 (1978), no. 116, 443--462.

\bibitem{Ingham}
Ingham, Albert Edward. 
Mean-Value Theorems in the Theory of the Riemann Zeta-Function. 
\emph{Proc. London Math. Soc.} (2) 27 (1927), no. 4, 273--300.

\bibitem{I-K04}
Iwaniec, Henryk; Kowalski, Emmanuel.
Analytic number theory. American Mathematical Society Colloquium Publications, 53. American Mathematical Society, Providence, RI, 2004. xii+615 pp. ISBN: 0-8218-3633-1

\bibitem{J-L21}
Jiang, Yujiao; L\"u, Guangshi.
The generalized Bourgain-Sarnak-Ziegler criterion and its application to additively twisted sums on $\rm GL_m$.
\emph{Sci. China Math.} 64 (2021), no. 10, 2207--2230.

\bibitem{J-L-W21}
Jiang, Yujiao; L\"u, Guangshi; Wang, Zhiwei.
Exponential sums with multiplicative coefficients without the Ramanujan conjecture.
\emph{Math. Ann.} 379 (2021), no. 1--2, 589--632.

\bibitem{KS00}
Keating, Jonathan P.;  Snaith, Nina Claire. 
Random matrix theory and $\zeta(1/2+it)$. 
\emph{Comm. Math. Phys}. 214 (2000), no. 1, 57--89. 

\bibitem{KS03}
Kim, Henry H. 
Functoriality for the exterior square of $\GL_4$ and the symmetric fourth of $GL_2$. 
With appendix 1 by Dinakar Ramakrishnan and appendix 2 by Kim and Peter Sarnak. 
\emph{J. Amer. Math. Soc.} 16 (2003), no. 1, 139--183.
 
\bibitem{Lu09}
L\"u, Guangshi.
On sums involving coefficients of automorphic $L$-functions.
\emph{Proc. Amer. Math. Soc.} 137 (2009), no. 9, 2879--2887.

\bibitem{LRS}
Luo, Wenzhi; Rudnick, Ze\'ev; Sarnak, Peter Clive.
On Selberg's eigenvalue conjecture.
\emph{Geom. Funct. Anal.} (1995), no. 2, 387--401.

\bibitem{Alex23}
Mangerel, Alexander P.
Divisor-bounded multiplicative functions in short intervals.
\emph{Res. Math. Sci.} 10 (2023), no. 1, Paper No. 12, 47 pp.

\bibitem{M-R15}
Matom\"aki, Kaisa; Radziwi\l\l, Maksym.
A note on the Liouville function in short intervals.
\emph{ArXiv preprints}
arXiv:1502.02374

\bibitem{M-R16}
Matom\"aki, Kaisa; Radziwi\l\l, Maksym.
Multiplicative functions in short intervals.
\emph{Ann. of Math.} (2) 183 (2016), no. 3, 1015--1056

\bibitem{M-R20}
Matom\"aki, Kaisa; Radziwi\l\l, Maksym.
Multiplicative functions in short intervals II.
\emph{ArXiv preprints}
arXiv:2007.04290

\bibitem{Nelson21}
Nelson, Paul D.
Bounds for standard $L$-functions.
\emph{ArXiv preprints}
arXiv:2109.15230

\bibitem{Ng21}
Ng, Nathan C.
The sixth moment of the Riemann zeta function and ternary additive divisor sums. 
\emph{Discrete Anal.} 2021, Paper No. 6, 60 pp.

\bibitem{NSW}
Ng, Nathan C.; Shen, Quanli; Wong, Peng-Jie.
The eighth moment of the Riemann zeta function
\emph{ArXiv preprints}
arXiv:2204.13891

\bibitem{Pal22}
Pal, Sampurna. 
Second moment of degree three $L$-functions. 
\emph{Int. Math. Res. Not. IMRN} 2025, no. 15, rnaf234.

\bibitem{Pi}
Pi, Qinghua. 
Fractional moments of automorphic $L$-functions on $\GL(m)$. 
\emph{Chinese Ann. Math. Ser. B} 32 (2011), no. 4, 631--642.

\bibitem{Ramachandra}
Ramachandra, Kanakanahalli. 
Some remarks on the mean value of the Riemann zeta function and other Dirichlet series. II. 
\emph{Hardy-Ramanujan J.} 3 (1980), 1--24. 

\bibitem{Shiu1980}
Shiu, Peter.
A Brun--Titchmarsh theorem for multiplicative functions.
\emph{J. Reine Angew. Math.} 313 (1980), 161--170.

\bibitem{Sound09moments}
Soundararajan, Kannan.
Moments of the Riemann zeta function. 
\emph{Ann. of Math.} (2) 170 (2009), no. 2, 981--993.
 
\bibitem{Sound10}
Soundararajan, Kannan.
Weak subconvexity for central values of $L$-functions.
\emph{Ann. of Math.} (2) 172 (2010), no. 2, 1469--1498.

\bibitem{S-T19}
Soundararajan, Kannan; Thorner, Jesse.
Weak subconvexity without a Ramanujan hypothesis. With an appendix by Farrell Brumley. 
\emph{Duke Math. J.} 168 (2019), no. 7, 1231--1268.

\bibitem{Sun24}
Sun, Yu-Chen
On divisor bounded multiplicative functions in short intervals.
\emph{ArXiv preprints}
arXiv:2401.08432 

\bibitem{TX16}
Tang, Hengcai; Xiao, Xuanxuan.
Integral moments of automorphic L-functions. 
\emph{Int. J. Number Theory} 12 (2016), no. 7, 1827--1843. 
 
\end{thebibliography}
\end{document}